\documentclass[journal]{IEEEtran}
\renewcommand{\baselinestretch}{1.0}
\usepackage{graphicx} 
\usepackage{mdwmath,mdwtab,url}
\usepackage{caption}
\usepackage{amsmath}
\usepackage{amssymb}
\usepackage{amssymb}
\usepackage{psfrag}
\usepackage{color}
\usepackage{cite}
\usepackage{epstopdf}

\ifCLASSINFOpdf
  
\else
  
\fi

\begin{document}
%
\title{Contour-FFT based Spectral Domain MBF Analysis of Large Printed Antenna Arrays}
%
%
%

\author{Shambhu~Nath~Jha,~\IEEEmembership{Member,~IEEE,}
        Christophe~Craeye,~\IEEEmembership{Senior Member,~IEEE}
\thanks{The authors are grateful to the R\'egion Wallone for financial support through the RADIANT and SUPLAN projects. The authors are with ICTEAM Institute, Universit\'e catholique de Louvain, Belgium ({shambhu.jha, christophe.craeye} @ uclouvain.be).}}

%
%

\markboth{Submitted to IEEE Transactions on Antennas and Propagation, on 25 March 2014}%
{}
%



\maketitle

\begin{abstract}
A fast spectral-domain method is proposed to evaluate the reaction terms between the Macro Basis Functions in regular and non-regular arrays made of identical printed antennas. The presented technique first exploits the filtering capabilities of the Macro Basis Functions in the spectral domain. The method is then strongly accelerated with the help of a newly formulated Fast Fourier Transform-based technique, which is applicable to a contour integration in the complex plane. We name the method as Contour-FFT or C-FFT. Besides an effective homogeneous medium term treated with multipoles, a computational complexity of order $N$$log_{2}N$ is achieved for the the tabulation of substrate-related reaction terms for any possible relative positions. The complexity of the proposed method is independent from the complexity of the elements. Numerical results obtained with the proposed method are compared with those from a pre-validated reference solution based on the traditional Macro Basis Functions technique; an excellent agreement is observed. 
\end{abstract}

\begin{IEEEkeywords}
Method of Moments, spectral domain method, macro basis function, contour-FFT, contour deformation, large irregular arrays, printed antenna, Green's function.
\end{IEEEkeywords}

%
\IEEEpeerreviewmaketitle

\section{Introduction}
\IEEEPARstart{T}{he Method-of-Moments} (MoM) \cite{Harng_book} technique is a very robust numerical method for the analysis of radiation and scattering problems involving large periodic or non-periodic structures: antenna arrays \cite{garcia}, reflectarrays \cite{reflect_book}, frequency selective surfaces \cite{Kasha}, and recently, metamaterials \cite{Radu_book}. The computation time needed to solve large electromagnetic problems by using a direct full-wave solution technique can be prohibitively high: the complexity of matrix filling increases as $O(A^p r B^2)$, where $A$ is the number of the antennas in the array, $B$ is the number of basis functions used to represent each element, and $r$ is the average complexity of computing one reaction integral. The value of $p$ is $1$ for regular arrays and $2$ for the case of irregular arrays. Besides, the computational complexity of solving the system of equations using a direct full-wave MoM technique scales as $O(BA)^3$. To solve large array problems, fast iterative techniques \cite{Janpugdee,itr1} have been proposed in the literature. Unfortunately, the number of required iterations is not known beforehand and the iterations need to be restarted everytime the array configuration or excitations change. Recently, various methods aim at decomposing the current on a given antenna into a limited set of current distributions covering the whole element, so as to dramatically reduce the number of unknowns at the element level. This domain-decomposition technique has led to the Macro Basis Functions (MBFs) method \cite{eucap2012:mosig}, the Characteristic Basis Functions Method (CBFM) \cite{eucap2012:mitra1}, the Synthetic Functions eXpansion (SFX) \cite{eucap2012:vecchi} and the Array Scanning Method-Macro Basis Function (ASM-MBF) \cite{CraeyeASM}. The relationship between this class of methods and iterative techniques has been studied in more detail in \cite{TAP2011Craeye}. These MBF-based methods reduce the computational complexity related to the solution time from $O(BA)^3$ to $O(MA)^3$, where $M$ $(M<<B)$ is the number of Macro Basis Functions used per element. However, a priori, using these fast methods, the computational complexity needed to fill the impedance matrix remains the same as that of a full-wave method. A method for spectral-domain interactions between analytically derived characteristic basis functions (CBFs) is proposed in \cite{ChMitrafast} to treat scattering problems from electrically large faceted bodies. In \cite{eucap2012:shambhuEucap}, an approach based on a Laurent Series model is proposed for the fast computation of MBF interactions in $1D$ irregular arrays of printed antennas. This work yields a computational complexity for the reduced-matrix filling that is independent from the element's complexity in the array, at the expense of some fixed amount of preparation time. Furthermore, the work in \cite{eucap2012:Ivica} presents an efficient way of computing MBF reaction terms and reduces them to equivalent moments. In \cite{mascant}, the CBFM is combined with Adaptive Cross Approximation algorithm to compute fast reduced interactions in electrically large arrays. The work in \cite {DavidInterp} proposes an interpolatory approach to compute the MBF interactions using a low-order harmonic-polynomial model obtained by considering a few physical transformations on the interaction function, such that very few calculations of the functions are needed. 

The present paper concerns a planar-spectral domain approach for computing MBF interactions. In \cite{jhaspectral}, we exploited the benefits of a planar spectral-domain approach and presented preliminary results related to the interactions between MBFs for identical elements in large irregular arrays of printed antennas. We further tested the application of the spectral-domain MBF method for the case of reflectarrays with non-identical elements \cite{jhaspectral_nonidentical}. In this paper, this class of methods is strongly accelerated using the so called Contour-FFT approach and is applied to regular and irregular arrays of identical elements. The use of the FFT itself for fast planar array analysis has been reported long time ago. For instance, the works \cite{FengFastFW,BingFastMethod} include the precorrected-FFT to speed up matrix-vector multiplications in iterative solutions. The authors of \cite{VechiAIM} proposed and demonstrated the use of the Adaptive Integral Method (AIM) fast factorization to accelerate the Synthetic Function eXpansion (SFX) domain decomposition by exploiting the convolutional nature of the Toeplitz kernel using a 3D FFT. This method is demonstrated to be efficient for volume and quasi-planar problems. Furthermore, the method proposed in \cite{NingFastMethod} combines the precorrected-FFT method and the discrete complex image method applicable to the Mixed Potential Integral Equation (MPIE) formulation. However, as far as the application of the 2D FFT is concerned, important accuracy problems arise from the treatment of surface-wave poles. For instance, when applying the complex image method as done in \cite{NingFastMethod}, this requires the accurate calculation of the contributions of the quasi-dynamic images, surface waves poles and complex images, leading to complicated formulations and possibly to approximated results if one misses some of those contributions. 

The problem has already been mentioned in $1990$ \cite{Catedra_FFT}, where one deals with the surface wave poles by using many integration points \cite{Catedra_FFT} around the poles, assuming that one has the knowledge of the exact location and of the number of poles, which is not really a trivial task.  The authors of \cite{ChMitrafast} suggest to use a fine spectral sampling near the singularities of the Green's function in spectral-domain. Another solution consists of neatly extracting the surface wave poles \cite{Pozar_Input},\cite{MosigPole}. A more general approach could exploit a contour deformation in complex wavenumber plane to avoid the surface wave poles that appear in the spectral-domain Green's function, as is done in many classical papers \cite{MosigContourDyadic}, \cite{Michalski}. Here, we will use a parabolic type of contour deformation \cite{contour_parabolic}. Unfortunately, with the application of complex-plane integration, as detailed in Section \ref{sec:Cfft}3, the integral kernel does no longer have a purely complex-exponential structure, i.e. the integral does no longer exactly correspond to a Fourier transform. This prevents the accurate acceleration of computations using the FFT algorithm. In this paper, we present a new formulation consisting of expanding the kernel obtained under contour deformation into a form suitable to FFT acceleration. In addition, the proposed Contour-FFT, as we name it, or C-FFT, is combined with a multipole-based MBFs interaction method\cite{ChFMM} dealing with the effective homogeneous medium part of the Green's function.  This multipole component is common to \cite{davidMP}. The implementation details of the multipoles-based MBF interactions are not reproduced in this paper; further information about them can be found in \cite{ChFMM} and \cite{davidMP}.

The C-FFT approach is used to very rapidly tabulate the substrate-related interactions between MBF versus relative positions, a goal that was also persued via a completely different method in \cite {DavidInterp} in the free-space case. This means that, once the table is constructed for a given substrate and a given antenna, the corresponding component of the reduced system of equation can be quasi-instantly constructed for any array configuration. This methodology offers very important advantages when mutual coupling needs to be included in the optimization of antenna arrays.

The remainder of this paper is organized as follows. Section \ref{sec:MBF} presents a brief reminder regarding the Macro Basis Function approach applied to antenna arrays printed on a substrate backed by a ground plane. In Section \ref{sec:fastspectral}, we first describe the spectral-domain approach for MBFs and we then present the C-FFT based fast spectral-domain method for Macro Basis Functions, applicable to large regular and irregular arrays of identical elements. In Section \ref{sec:Results}, validation results are shown and the computation times and computational complexity of the proposed method are also reported. Finally, conclusions are drawn in Section \ref{sec:conclusion}.

\section{MBF Approach for Printed Antennas}
\label{sec:MBF}
Following the MoM technique, the final system of equations for an antenna array is given in the form 
\begin{equation}
	\underline{\underline{Z}} \hspace{0.1cm} \underline{I}=\underline{V} 
\end{equation}
where $\underline{\underline{Z}}$ is the MoM interaction matrix of size $U \times U$ ($U$ being the total number of unknowns) consisting of the interactions between all elementary basis functions in the array. $\underline{I}$ is an unknown coefficient vector with size $U \times 1$ and $\underline{V}$ is an excitation vector of size $U \times 1$. Our full-wave Method of Moment code is based on the Electric-Field Integral Equation (EFIE) \cite{MichalskiEFIE}, \cite{mosigEFIE} and the MoM impedance matrix is given as:
\begin{equation}
\underline{\underline{Z}}=\underline{\underline{Z}}^l+\underline{\underline{Z}}^h
\end{equation}
where, $\underline{\underline{Z}}^l$ is the contribution from layered-medium Green's function, excluding the contribution from an effective homogeneous medium \cite{ValerioExt} and its image \cite{PozarImage}. $\underline{\underline{Z}}^h$ is the matrix for the effective homogeneous medium and its image, which can be directly estimated through space-domain convolution.   

The MBF method reduces the size of the MoM system of equations by replacing the original set of elementary basis functions by a set of aggregate basis functions. ``Macro Basis Functions" \cite{eucap2012:mosig} are obtained through the solution of smaller problems and can be generated in various ways \cite{eucap2012:mosig} - \cite{CraeyeASM}. In this work, the MBFs based on ``primary-and-secondaries" approach proposed in \cite{eucap2012:mitra1} is used. 
Let us denote by $Z_{mn}$ the block of $Z$ corresponding to basis functions on antenna $n$ and testing functions on antenna $m$. If $Q$ is the matrix whose columns describe the MBFs, then the corresponding block of
the reduced matrix can be obtained as follows:
\begin{equation}
\label{eq:MBF}
Z'_{mn}=Q^HZ_{mn}Q 
\end{equation}
Let us denote by $V_m$ the segment of $V$ corresponding to testing functions on antenna $m$. Then, the corresponding segment of the reduced excitation vector reads:
\begin{equation}
\label{eq:voltg}
V'_m=Q^H V_m 
\end{equation}
where superscript ($^H$) in (\ref{eq:MBF}) and (\ref{eq:voltg}) represents the Hermitian transposed.

With the MBF-based method, the original matrix of size $BA \times BA$ is reduced into $MA \times MA$. Usually, $M$ is very small compared to $B$. Hence, the MBF method often allows solving the large system of equations by direct inversion without resorting to iterative techniques. In this work, $9$ MBFs (1 primary and 8 secondaries) have been used. Our goal is to accelerate the evaluation of $Z'_{mn}$ for the case of antennas printed on a layered medium.

\section{Spectral-Domain Approach}
\label{sec:fastspectral}
We first recall the spectral-domain Method of Moments (MoM) approach \cite{Itoh_book,Spectral} for interactions between elementary basis functions. The underlying formulation is then extended to compute the interactions between MBFs. We demonstrate that by applying the spectral-domain approach to MBFs, one can directly obtain the reduced matrix for the MBF-based method, with a saving by a large factor in the spectral integration time by exploiting the filtering properties of MBFs. We then explain how the C-FFT can be applied to expedite the evaluation of spectral-domain MBF interactions such that an $Nlog_2N$ complexity is achieved, where, $N$ is the number of points used in the FFT. The details on $N$ is provided in Section \ref{sec:comp}.
\begin{figure}[!t]
\begin{center}
\includegraphics[width=0.49\textwidth]{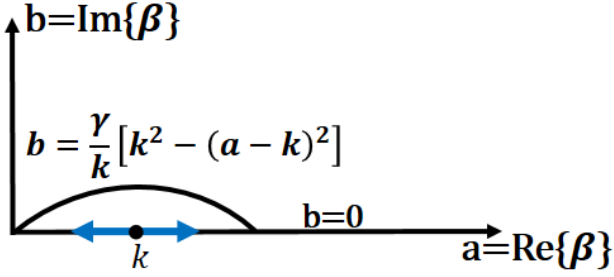}
\caption{1D Contour deformation technqiue for complex plane integration. The deformed contour is of parabolic type with a factor $\gamma$, a ratio corresponding to the contour height. The arrow in the X-axis of the plot represents the possible location of surface wave poles on the real axis.}
\label{fig:contour_1D}
\end{center}
\end{figure}
\subsection{Spectral-Domain based MBF Interactions}
\label{sec:FFT-spec}
The spectral-domain formulation for interaction between elementary basis functions is given as \cite{Spectral}
\begin{equation}
\label{eq:spectral}
Z_{tb}=\frac{1}{4\pi^2}\int \int \tilde{H} e^{-i(k_x\Delta x+k_y\Delta y)} dk_x dk_y \\
\end{equation} 
\begin{equation}
\label{eq:spectralkernel}
\tilde{H}=\tilde{\underline{\underline{F_t}}}\cdot \tilde{\underline{\underline{G}}} \cdot \tilde{\underline{\underline{F_b}}}^*
\end{equation} 
$\tilde{\underline{\underline{F_t}}}$ and $\tilde{\underline{\underline{F_b}}}$ are the Fourier transforms of testing and basis functions, respectively. The Fourier transforms of rooftop basis functions, used below, are available in analytical form. The terms $\Delta x$ and $\Delta y$ give the distance between the source and observation domains, taken here as the antennas on which basis and testing functions are respectively located. The centers of those domains serve as local origins for the calculation of $\tilde{\underline{\underline{F_t}}}$ and $\tilde{\underline{\underline{F_b}}}$. $\tilde{\underline{\underline{G}}}$ is the Dyadic form of the layered medium Green's function \cite{Pozar_Input, MosigContourDyadic, eucap2010:dyadicgreen}, from which an effective homogeneous medium \cite{ValerioExt} contribution and its image through the ground plane \cite{PozarImage} have been extracted. The formulation (\ref{eq:spectral}) also holds when the elementary basis functions are replaced by MBFs; this is then a numerical extension of the reaction integrals between analytically known CBFs\cite{ChMitrafast}.

\begin{figure}[!t]
\psfragscanon
\begin{center}
\psfrag{y}[][][0.8]{Real($k_y$)}
\psfrag{x}[][][0.8]{Real($k_x$)}
 \includegraphics[height=5.8cm, width=0.47\textwidth]{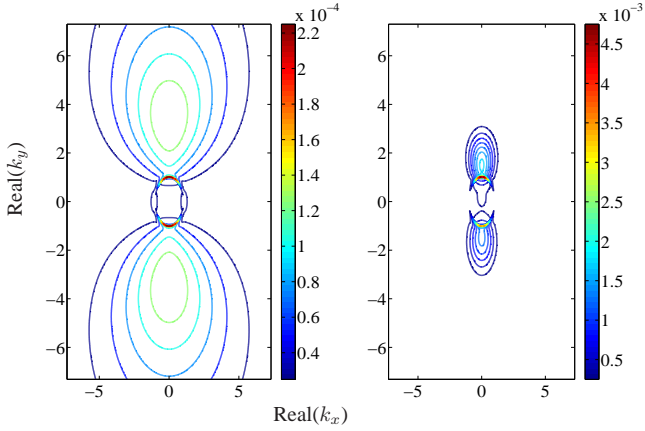}
\end{center}
\caption{Spectral Integrands. Left: integrand with elementary basis function, Right: integrand with MBFs.}
    \label{fig:integrands_spectral}
\end{figure}

For the sake of computational simplicity, the above integration (\ref{eq:spectral}) is carried out in polar coordinates, through the following change of variables: $k_x=\beta \cos\alpha$, $k_y=\beta \sin\alpha$, $\rho=\sqrt{{\Delta x}^2+{\Delta y}^2}$ and $\phi=\arctan({\Delta y}/{\Delta x})$. This leads to
\begin{equation}
\label{eq:polarcord}
Z_{tb}=\frac{1}{4\pi^2}\int^{2\pi}_{\alpha=0}\int^{\infty}_{\beta=0} \tilde{H} e^{-i\beta \rho \cos(\alpha-\phi)} \beta d\beta d\alpha \\
\end{equation} 
As pointed out in the introduction, a parabolic type of contour deformation \cite{contour_parabolic} as shown in Fig. \ref{fig:contour_1D} is applied along the $\beta$ coordinate to avoid poles that may appear in the Green's function along the real axis. The same deformation remains valid for all values of the $\alpha$ angular coordinate. The work in \cite{Radio_inhomogeneous} also considers a contour deformation for 3D problems in the free-space case, with integration along $\theta=arc\sin(\beta/k_0)$ instead of $\beta$, where $k_0$ is a wavenumber in free-space. The formulation (\ref{eq:polarcord}) is updated below to include the effect of contour deformation in complex $\beta$ plane with  $\beta=\beta_c=\beta_R+i\beta_I$, while explicitly integrating along $\beta_R$. The resulting MoM matrix entry considering integration along the contoured path becomes: 
\begin{equation}
\label{eq:spectral_contour}
Z_{tb}=\frac{1}{4\pi^2}\int_{\alpha}\int_{\beta}\left(1+i\frac{d\beta_I}{d\beta_R}\right) \tilde{H} e^{-i\beta_c \rho \cos(\alpha-\phi)} \beta_c  d\beta_R  d\alpha \\
\end{equation} 
The formulation of (\ref{eq:spectral_contour}) can be extended to compute the interactions between MBFs for any relative position of the antennas in the array as:
\begin{equation}
\label{eq:MBF_mul}
Z'_{mn}=\frac{1}{4\pi^2} \int \int \left(1+i\frac{d\beta_I}{d\beta_R}\right) \tilde{H}_M e^{-i\beta_c \rho \cos(\alpha-\phi)} \beta_c  d\beta_R  d\alpha
\end{equation}  
\begin{equation}
\label{eq:MBF_mul22}
\tilde{H}_M=\tilde{\underline{\underline{M_t}}}\cdot \tilde{\underline{\underline{G}}}\cdot \tilde{\underline{\underline{M_b}}}^* 
\end{equation}  
where $\tilde{\underline{\underline{M_t}}}$ and $\tilde{\underline{\underline{M_b}}}$ are the Fourier transforms of macro-testing and macro-basis functions, respectively. Let us consider $q_j$ as the coefficients multiplying the $j^{th}$ elementary basis function in the MBF description and $\tilde{\underline{\underline{F_j}}}$ as the Fourier transform of that elementary basis functions. The Fourier transform of MBFs is then given as:
\begin{equation}
\label{eq:FTMBF}
\tilde{\underline{\underline{M}}}=\sum_j {q_j}^* \tilde{\underline{\underline{F_j}}}
\end{equation}
As can be demonstrated from (\ref{eq:MBF_mul}) and (\ref{eq:MBF_mul22}), the reaction terms between MBFs can be evaluated by a simple multiplication of the Fourier transforms (FT) of Macro Basis Functions (MBFs) and Macro Testing Functions (MTFs), and the Fourier transforms of the Dyadic form of the Green's function for the layered medium. It should be well noted here that the contribution from an effective homogeneous medium part of the Green's function and its image through the ground plane are treated with the multipole approach proposed in \cite{ChFMM}. In the proposed method, we exploit the benefit that, as the MBFs have a wider support in space domain, their Fourier transforms are narrow functions. It leads to a rapidly decaying spectral integrand, hence allowing a limitation of the integration domain. The filtering effect of the Fourier transform of MBFs on the spectral integrands for a microstrip type of antenna at $24.125$ GHz is illustrated in Fig. \ref{fig:integrands_spectral}. The length of the antenna used for simulation is $3.82$ mm and the width is $6.5$ mm. To properly match the antenna, a quarter-wave matching transmission line is used. The length and width of the quarter-wave transmission line are $2.4$ mm and $0.64$ mm, respectively. The antenna is printed on a RT Duroid substrate of relative permitivity $\epsilon_r=2.2$ and thickness $h=0.381$ mm. The antenna is meshed with $243$ rooftop basis functions. As can be seen, the spectral integrands in the case of MBF-based interactions are very narrow functions, as compared to the integrand obtained for elementary basis functions. Let us denote by $k_d$ the wavenumber in the dielectric material. With the help of the filtering properties of MBFs in spectral-domain, a spectral integration limit of $5k_d$ has been found sufficient to obtain accurate results for interactions between MBFs. In this way, it is possible to save a large factor in the filling time of the reduced MoM matrix. Following the Galerkin testing procedure, the sets of MBFs and MTFs are the same. 

\subsection{C-FFT based Formulation for Fast MBF Interactions}
\label{sec:Cfft}
\subsubsection{Tabulation}
Although the formulation of (\ref{eq:MBF_mul}) is quick, for the case of large irregular arrays, one needs to compute order of $A^2$ interactions, where $A$ is the number of antennas. Despite the filtering effect reported in the previous section, the computation time is still prohibitive. For large irregular arrays, it would be necessary to evaluate (\ref{eq:MBF_mul}) for every relative position between elements in the array. We would like to obtain a method which allows us to pre-compute the MBF interactions for any arbitrary array spacing without explicitly repeating everytime a similar evaluation. In Section \ref{sec:fftmethod}, we develop a method which allows the tabulation of the MBF interactions in space-domain over a very fine grid of relative distances $(\Delta x,\Delta y)$ between antennas. In this respect, the FFT will be extremely useful, provided that appropriate transformations are made to accomodate the effects of contour deformation. Then, to fill the reduced MoM impedance matrix for a given array configuration, the interactions between MBFs at arbitrary spacings can be computed extremely fast through a low-order interpolation in the table. 
\subsubsection{Polar-to-Rectangular Coordinates Transformation}
To exploit the speed-up factor offered by the FFT algorithm, a necessary condition is to come back to a regular sampling in the rectangular $k_x-k_y$ coordinates, still with proper treatment of the surface-wave poles through contour deformation. We propose a formulation which allows us to keep integrating along the real values of $k_x-k_y$ in rectangular coordinates and which includes the contour deformation. As we will see, this transformation in itself, is not sufficient for a direct application of the FFT; this problem will be solved in the next section. 

As results from the Section \ref{sec:FFT-spec}, the 1D contour deformation technique is first extended to cover the 2D ($k_x$, $k_y$) plane. As shown below, the 2D contour deformation can be equally appplied in the $k_x-k_y$ domain by first finding a suitable complex-plane extension for $k_x$ and $k_y$. 

Let us denote the corresponding complex transverse wavenumbers as $k_x=k_{xR}+ik_{xI}$ and $k_y=k_{yR}+ik_{yI}$. They have to satisfy 
\begin{equation}
\label{eq:complexbeta}
({k_{xR}+ik_{xI}})^{2}+({k_{yR}+ik_{yI}})^{2}=(\beta_R+i\beta_I)^2=k^2_0-k^2_{z}
\end{equation}
where $k_{z}$ is the wavenumber along $z$ used in the Green's function.
Expanding (\ref{eq:complexbeta}) into its real and imaginary parts leads to the following relations:
\begin{equation}
\label{eq:rel1}
\beta^2_R-\beta^2_I={k^{2}_{xR}}-{k^{2}_{xI}}+{k^{2}_{yR}}-{k^{2}_{yI}}
\end{equation}
\begin{equation}
\label{eq:rel2}
\beta_R\beta_I=k_{xR}k_{xI}+k_{yR}k_{yI}
\end{equation}
By definition, we also have 
\begin{equation}
\label{eq:beta}
{k^2_{xR}+k^2_{yR}}=\beta_R^2
\end{equation}
Equation (\ref{eq:rel2}) can be rewritten as:
\begin{equation}
\label{eq:rel2new}
\beta^2_R \frac{\beta_I}{\beta_R}= k^2_{xR} \frac{k_{xI}}{k_{xR}} + k^2_{yR} \frac{k_{yI}}{k_{yR}} 
\end{equation}
Subtracting (\ref{eq:beta}) from (\ref{eq:rel1}) provides, after a simple transformation
\begin{equation}
\label{eq:rel1new}
\beta^2_R {\left(\frac{\beta_I}{\beta_R}\right)}^2= {\left(\frac{k_{xI}}{k_{xR}}\right)}^2 k^2_{xR} + {\left(\frac{k_{yI}}{k_{yR}}\right)}^2  k^2_{yR}
\end{equation}
By inspection of (\ref{eq:rel2new}) and (\ref{eq:rel1new}), it becomes clear that the following rule for the analytical extension satisfies (\ref{eq:complexbeta}):
\begin{equation}
\label{eq:ratioNew}
\frac{k_{xI}}{k_{xR}}=\frac{k_{yI}}{k_{yR}}=\frac{\beta_I}{\beta_R} \triangleq \gamma
\end{equation}
Writing (\ref{eq:MBF_mul}) back into the rectangular $k_{xR}-k_{yR}$ domain leads to:
\begin{equation}
\label{eq:FFT}
Z'_{mn}=\frac{1}{4\pi^2}\int\int \tilde{I}_M e^{-i(k_{xR}\Delta x+k_{yR}\Delta y)} E_I dk_{xR} dk_{yR} \\
\end{equation} 
\begin{equation}
E_I = e^{\gamma(k_{xR} \Delta x+k_{yR} \Delta y)}
\end{equation} 
\begin{equation}
\tilde{I}_M=\left(1+i\frac{d\beta_I}{d\beta_R}\right) (1+i\gamma) \tilde{H}_M
\end{equation} 
where the factor $(1+i\gamma)$ results from the following identity: $\beta_c d\beta_R d\alpha=(1+i\gamma)dk_{xR}dk_{yR}$. 
Formulation (\ref{eq:FFT}) allows us to transform the integration along the complex polar coordinates into an integration along real $k_x-k_y$ values; non-zero values of $\gamma=\beta_I/\beta_R$ and $d\beta_I/d\beta_R$ account for contour deformation in the complex plane. 
\subsubsection{C-FFT Based Method}
\label{sec:fftmethod}
If one has a rapid look at the structure of (\ref{eq:FFT}), it can be seen that it almost has a Fourier transform structure, the difference being due to the real exponential $E_I$ arising from the contour deformation. Hence, in this form, the FFT algorithm cannot be directly applied. A solution to this problem comes with the Taylor's series expansion of $E_I$. 

From now on, for simplicity of notation, subscript $_R$ in $k_{xR}$ and $k_{yR}$ will be omitted, such that a spectral pair ($k_x$,$k_y$) will always denote the real values. For a spectral pair ($k_x$,$k_y$) and spatial distances ($\Delta x$,$\Delta y$), the Taylor's series expansion of factor $E_I$ reads,
\begin{equation}
\label{eq:taylors}
E_I=1+\gamma(k_{x}\Delta x+k_{y}\Delta y)+\frac{\gamma^2}{2}(k_{x}\Delta x+k_{y}\Delta y)^2+..
\end{equation} 
It should be noted that the convergence of the Taylor's series (\ref{eq:taylors}) will be dictated by $\gamma k_0 d_{max} << 1$, where $d_{max}$ is the maximum spacing of the array and $\gamma=\beta_I/\beta_R$ measures the height of the deformed contour. Hence, the convergence will be faster for lower contours and smaller arrays in terms of wavelengths. As will be explained in Sections \ref{Taylors} and \ref{sec:reduced}, in view of the low contours that can be afforded, it will be possible to analyse large arrays using very few terms in the series.

Inserting (\ref{eq:taylors}) into (\ref{eq:FFT}), while taking separately the terms of the Taylor's series, and moving outside the integral all the factors depending on $\Delta x$ and $\Delta y$ leads to:
\begin{eqnarray}
\label{eq:fftexp}
Z'_{mn}=\frac{1}{4\pi^2}\sum_p\sum_q \tilde{I}_M e^{-i(k_{x}\Delta x_p+k_{y}\Delta y_q)} \Delta k_{x}\Delta k_{y} \nonumber \\
+\Delta x \frac{1}{4\pi^2}\sum_p\sum_q \gamma k_{x} \tilde{I}_M e^{-i(k_{x}\Delta x_p+k_{y}\Delta y_q)}\Delta k_{x} \Delta k_{y} \nonumber \\
+\Delta y \frac{1}{4\pi^2}\sum_p\sum_q \gamma k_{y} \tilde{I}_M e^{-i(k_{x}\Delta x_p+k_{y}\Delta y_q)}\Delta k_{x} \Delta k_{y}+..
\end{eqnarray}
where integrals have been approximated by summations, according to the basic rectangles rule.

Now we can see that the summations exactly have a DFT structure. It then follows that the FFT algorithm can be directly applied for each term of the Taylor's series expansion. Stemming from (\ref{eq:taylors}), the space domain parameters $\Delta x$ and $\Delta y$ are factored out, which allows us to compute the interactions at any relevant array spacing, once we have the result of the FFT operations. 

Writing (\ref{eq:fftexp}) into a standard FFT form considering a Taylor's series up to third order leads to the following relation:
\begin{eqnarray}
\label{eq:fftarry}
Z'_{mn} (k_1,k_2) \approx K'[FFT2(\tilde{I}_{Mf})+\Delta xFFT2(\gamma k_{x}\tilde{I}_{Mf})\nonumber \\
+\Delta yFFT2(\gamma k_{y}\tilde{I}_{Mf})+\frac{1}{2}(\Delta x)^2FFT2(\gamma^2 k^2_{x} \tilde{I}_{Mf})\nonumber \\
+(\Delta x \Delta y)FFT2(\gamma^2 k_{x}k_{y} \tilde{I}_{Mf})\nonumber \\
+\frac{1}{2}(\Delta y)^2FFT2(\gamma^2 k^2_{y} \tilde{I}_{Mf})+\frac{1}{6}(\Delta x)^3 FFT2(\gamma^3 k^3_{x} \tilde{I}_{Mf})\nonumber \\
+\frac{3}{6}(\Delta x)^2(\Delta y) FFT2(\gamma^3 k^2_{x} k_y \tilde{I}_{Mf})\nonumber \\
+\frac{3}{6}(\Delta x)(\Delta y)^2 FFT2(\gamma^3 k_x k^2_{y} \tilde{I}_{Mf})\nonumber \\
+\frac{1}{6}(\Delta y)^3 FFT2(\gamma^3 k^3_{y} \tilde{I}_{Mf})]
\end{eqnarray}
 where $k_1$ and $k_2$ are the space-domain indices.

The expression for $K'$ follows from the change of variables necessary to conform the DFT to the FFT structure. The details regarding the variables changes are provided in Appendix I. $\tilde{I}_{Mf}$ is a modified form of $\tilde{I}_{M}$; the details of this modification are given in the same Appendix. As demonstrated in Section \ref{sec:reduced}, by considering second or third order terms in the expansion, a very good accuracy level can be obtained for the MBFs interactions, as compared to those provided by the traditional MBF method. 

The formulation (\ref{eq:fftarry}) is directly applicable to the case of irregular arrays of identical elements. Once the values of $Z'_{mn}$ are computed using the FFT, one can use an interpolation technique to obtain the reduced MBF interactions for arbitrary values of ($\Delta x$,$\Delta y$). The interpolation can be limited to a low order of the tables produced by the FFT over a sufficiently fine grids. Therefore, we make use of the Whittaker-Shannon interpolation formulation with zero-padding, leading to the double of the originally considered FFT size. Given this preliminary refinement of the tables generated with the C-FFT, when filling the reduced MoM matrix, it is sufficient to search in the tables, using a second order interpolation method.

\section{Results and Discussions}
\label{sec:Results}
Simulations have been carried out for arrays made of the microstrip antenna described before in Section \ref{sec:FFT-spec} at the frequency of $24.125$ GHz. At first, in Section \ref{sec:normal}, the results obtained with the code implemented with traditional MBFs is validated for the case of irregular arrays of $25$ elements, by comparison with the results obtained with direct inversion and with the commercial software IE3D \cite{IE3D}. Then, in Section \ref{Taylors}, the relationship between the contour heights and the Taylor's series orders is briefly illustrated in Table \ref{tab:contour}. In Section \ref{sec:reduced}, a reduced interaction $Z'_{mn}$ between two primary MBFs computed with the proposed C-FFT based MBF method (\ref{eq:fftarry}) are presented and compared with the results obtained with the traditional MBF approach.  In Section \ref{sec:irrgarrys}, we discuss the results in terms of port currents and radiation patterns for the case of irregular arrays made of $100$ identical microstrip patches. Further results regarding port currents are provided for the case of large regular arrays of size $25$ $\times$ $25$ in Section \ref{sec:regarras}. The methods are computationally benchmarked in both cases and simulation times are reported in Table \ref{tab:computation}. The simulations are carried out using a PC with $16$ GB RAM and Processor Intel (R) Core(TM) $i7-3770$ CPU @$3.4$ GHz. The errors shown in the results are given in dB and are defined as $10\log_{10}$$\mid$$I_a$$-$$I_r$$\mid^2$ where, $I_a$ and $I_r$ are the actual (C-FFT) and reference (Traditional MBF) values, respectively. The results are normalized w.r.t. the maximum value of $\left|I_r\right|$.
\subsection{Validation of the traditional MBF Method Implementation}
\label{sec:normal}
\begin{figure} 
\includegraphics[height=6.2cm, width=0.4\textwidth]{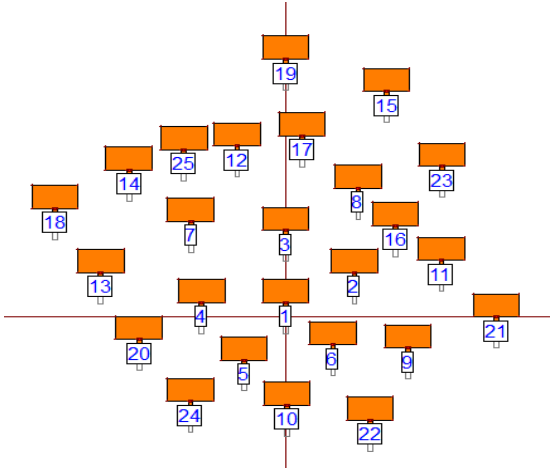}
\caption{An irregular array with 25 elements, plot generated with IE3D.}
\label{fig:Irrg25}   
\end{figure}

\begin{figure}[!t]
\psfragscanon
\begin{center}
\psfrag{a}[][][0.7]{\hspace{11em}{Port currents with IE3D}}
\psfrag{b}[][][0.7]{\hspace{16em}{Port currents with Fullwave Method}}
\psfrag{c}[][][0.7]{\hspace{14.5em}{Port currents with MBF Method}}
\psfrag{y}[][][0.8]{Port currents (dBA)}
\psfrag{x}[][][0.8]{Antenna Index in an Irregular Array with 25 Elements}
 \includegraphics[height=6.2cm, width=0.44\textwidth]{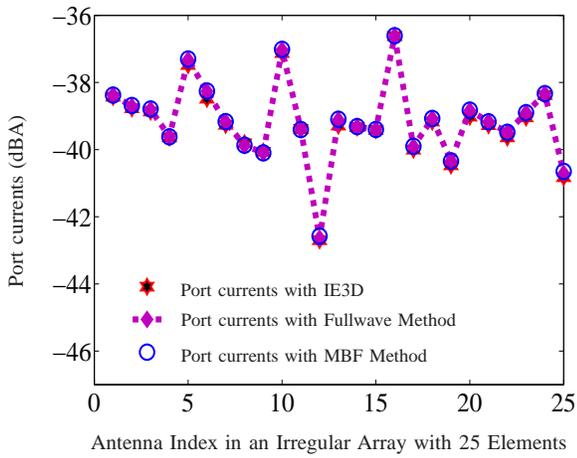}
\end{center}
\caption{Port currents in a 25 element irregular array when all the elements are excited.}
\label{fig:IE3D_val_port}
\end{figure}

\begin{figure}[!t]
\begin{center}
\psfrag{a}[][][0.7]{\hspace{9.5em}{Pattern with IE3D}}
\psfrag{b}[][][0.7]{\hspace{14.5em}{Pattern with Fullwave Method}}
\psfrag{c}[][][0.7]{\hspace{13em}{Pattern with MBF Method}}
\psfrag{y}[][][0.8]{Directivity (dBi)}
\psfrag{x}[][][0.8]{Angle from Broadside (deg)}
\includegraphics[height=5.8cm, width=0.44\textwidth]{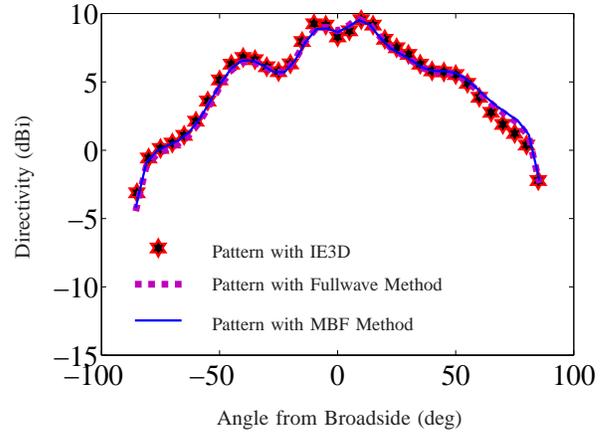}
\end{center}
\caption{Embedded element pattern when a single antenna in the array is excited.}
\label{fig:IE3D_val_pat}
\end{figure}

In this section, the implementation of the traditional MBF method based on the ``primary-and-secondaries" approach is first validated with the full-wave results and the results obtained with the commerical software IE3D\cite{IE3D}. For this purpose, an irregular array with $25$ elements located very close to each other, as shown in Fig. \ref{fig:Irrg25}, is considered. The port currents obtained when all the elements are excited in phase are shown in Fig. \ref{fig:IE3D_val_port}. As can be seen, this particular example exhibits very high mutual coupling effect as the variations in the port currents can be as high as $6$ dB. Moreover, an E-plane embedded element pattern is given in Fig. \ref{fig:IE3D_val_pat} for the case where antenna $1$ in the middle of the array is excited and the rest is terminated with $50$ $\Omega$ loads. These results validate the engineering accuracy of the implemented Method of Moments and MBF codes. Hence, the results obtained with this traditional MBF method code will be used as a reference solution to benchmark the performance of the C-FFT based MBF methods.

\subsection{On the Contour height and Taylor's series order}
\label{Taylors}
The evolution of the error in the final solution for various values of $\gamma$ [see equation (\ref{eq:ratioNew})] and for various orders of the Taylor's series terms [see equation (\ref{eq:fftarry})] is summarized in Table \ref{tab:contour}. These errors correspond to the computation of interactions between two primary MBFs on the antennas placed diagonally at the distance of $d=\sqrt{\Delta x^2+\Delta y^2}= 21.2\lambda_0$, where, $\lambda_0$ is the free-space wavelength. An FFT of size $2048$ $\times$ $2048$ is used for this error calculation. The errors are given in dB and are defined as before. There is a trade-off between the contour heights and the required number of terms needed. However, too high (e.g. $\gamma=1/50$) and too low (e.g. $\gamma=1/400$) values of $\gamma$ should be avoided. For instance, the higher the contour is made, the more Taylor's series terms are needed for better accuracy. Too low contours would make the surface-wave poles visible. In practice, any value of $\gamma$ between $1/100$ and $1/200$ gives accurate results by considering $2$nd or $3$rd order of Taylor's series terms. In the simulations described further, a $\gamma$ value of $1/130$ has been used. 
\begin{table}[ht]
\caption{Error level in dB versus contour height and order of Taylor's series.}  
\centering  
\begin{tabular}{c c c c c c} 
\hline 
\hline 
\emph{Order}/$\gamma$ &1/50 &1/100 &1/130 &1/200 &1/400\\[0.5ex]
\hline
0 & -0.3& -2.1 & -3.5 & -6.8 & --7.9\\
1 & -1.5& -6.1 & -8.6 & -11.9 & -6.6\\
2 & -4.0& -12.8 & -16.8& -19.4 & -7.3\\
3 & --7.7& -20.4 & -24.4& -22.0 & -7.3\\[1ex]
\hline
\end{tabular}
\label{tab:contour}
\end{table}
\subsection{Computation of the Reduced Matrix}
\label{sec:reduced}
\begin{figure}[!t]
\psfragscanon
\begin{center}
\psfrag{a}[][][0.7]{\hspace{10.75em}{Traditional MBF Method}}
\psfrag{b}[][][0.7]{\hspace{9.75em}{Zero Order Term only}}
\psfrag{c}[][][0.7]{\hspace{9.75em}{Upto 1st Order Terms}}
\psfrag{d}[][][0.7]{\hspace{10.0em}{Upto 2nd Order Terms}}
\psfrag{e}[][][0.7]{\hspace{10.0em}{Upto 3rd Order Terms}}
\psfrag{y}[][][0.8]{Normalized Interactions in dB}
\psfrag{x}[][][0.8]{Relative distance in terms of $\lambda_0$}
\includegraphics[height=6.4cm, width=0.49\textwidth]{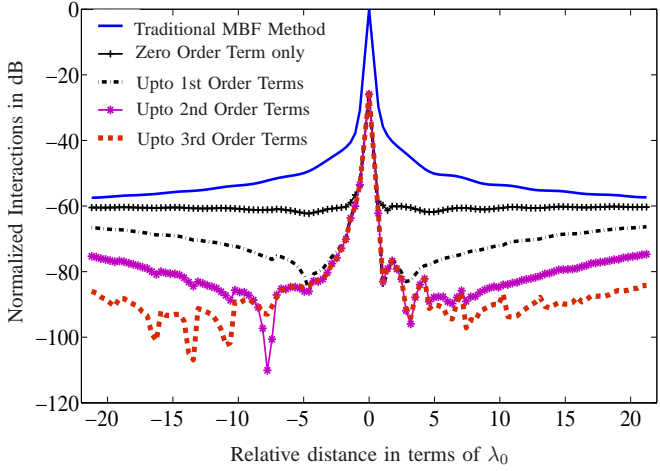}
\end{center}
\caption{Accuracy of the C-FFT based method in relation to the number of Taylor's series terms used. The results are computed with $4086 \times 4086$ FFT for better representation of the evolution of accuracy.}
    \label{fig:evolution_error_high}
\end{figure}

\begin{figure}[!t]
\psfragscanon 
\begin{center}
\psfrag{a}[][][0.7]{\hspace{10.75em}{Traditional MBF Method}}
\psfrag{b}[][][0.7]{\hspace{9.75em}{Zero Order Term only}}
\psfrag{c}[][][0.7]{\hspace{9.75em}{Upto 1st Order Terms}}
\psfrag{d}[][][0.7]{\hspace{10.0em}{Upto 2nd Order Terms}}
\psfrag{e}[][][0.7]{\hspace{10.0em}{Upto 3rd Order Terms}}
\psfrag{y}[][][0.8]{Normalized Interactions in dB}
\psfrag{x}[][][0.8]{Relative distance in terms of $\lambda_0$}
\includegraphics[height=6.4cm, width=0.49\textwidth]{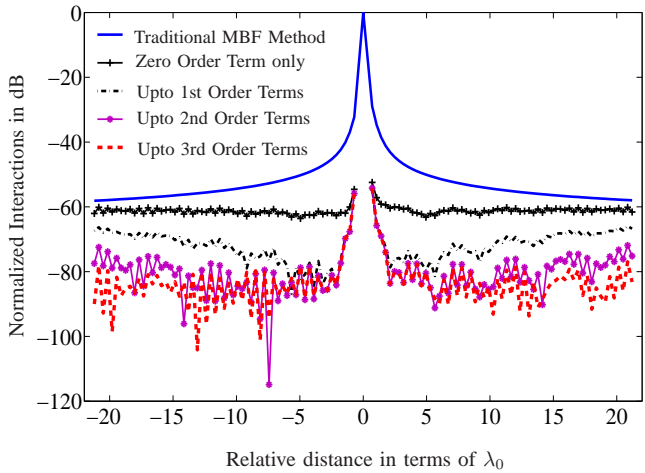}
\end{center}
\caption{Accuracy of the C-FFT based method in relation to the number of Taylor's series terms used. The results are full interactions: the layered medium part computed with FFT of size $2048 \times 2048$ and combined with the contribution from the effective homogeneous medium based on the MBF method for the multipole approach.}
    \label{fig:evolution_error}    
\end{figure}

In this section, we present the reduced matrix computed with the C-FFT based MBF method and we compare the results with those obtained with the traditional MBF method. The results are generated for various relative distances between the two elements placed diagonally. The relative distances in the following plots are given in terms of $\sqrt{\Delta x^2+\Delta y^2}$. $\Delta x$ and $\Delta y$ are the relative distances between the two elements in X and Y directions, respectively. The comparisons are given in Fig. \ref{fig:evolution_error_high}, which shows the evolution of the accuracy in the computation of the MBF interactions using the C-FFT based approach for different numbers of terms in the Taylor's series. The results are given in normalized value in dB scale. In this plot, an FFT of size $4096 \times 4096$ is deliberately used to obtain a better clarity in the evolution of errors. Similar plots are given in Fig. \ref{fig:evolution_error}, but for an FFT of size $2048 \times 2048$, on which all the remaining results are based on. The results in Fig. \ref{fig:evolution_error} also include the contributions from the effective homogeneous medium and its image, evaluated efficiently with the multipole approach applied to MBFs \cite{ChFMM}. The errors shown in the plots are given in dB and are defined as before. The errors in the computation of the results are below $-30$ dB. As can be seen from Figs. \ref{fig:evolution_error_high} and \ref{fig:evolution_error}, the accuracy in the result clearly increases as the number of terms in the Taylor's series is increased. It can be observed that just by considering up to second or third order terms, one can obtain very accurate results.  

\subsection{Analysis of Irregular Array of Microstrip Antennas}
\label{sec:irrgarrys}
\begin{figure}[!t]
\psfragscanon
\begin{center}
\psfrag{y}[][][0.8]{Y-distance in terms of $\lambda_0$}
\psfrag{x}[][][0.8]{X-distance in terms of $\lambda_0$}
 \includegraphics[height=5.3cm, width=0.44\textwidth]{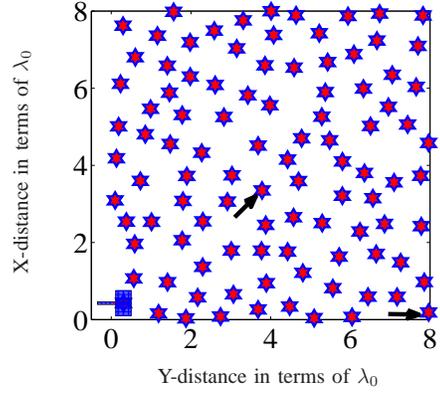}
\end{center}
\caption{Irregular arrays with 100 randomly placed elements.}
    \label{fig:array_irrg}  
\end{figure}

For this case, a total number of $100$ elements is taken, with irregular spacings, as shown in Fig. \ref{fig:array_irrg}. The elements are spread over a square area of size $8 \lambda_0$ $\times$ $8 \lambda_0$. The spacing between the elements varies from $0.0447 \lambda_0$ for the nearest elements to $7.86 \lambda_0$ for the farthest elements in the X-direction and from $0.1279 \lambda_0$ to $7.2525 \lambda_0$ in the Y-direction. The total number of unknowns is $24300$. 
\begin{figure}[!t]
\psfragscanon
\begin{center}
\psfrag{a}[][][0.7]{\hspace{11.5em}{Traditional MBF Method}}
\psfrag{b}[][][0.7]{\hspace{12.25em}{Error with zero order term}}
\psfrag{c}[][][0.7]{\hspace{11.75em}{Error with 1st order term}}
\psfrag{d}[][][0.7]{\hspace{12.0em}{Error with 2nd order term}}
\psfrag{e}[][][0.7]{\hspace{12.0em}{Error with 3rd order term}}
\psfrag{y}[][][0.8]{Normalized port currents (dB)}
\psfrag{x}[][][0.8]{Antenna Index in an Irregular Array with 100 Elements}
 \includegraphics[height=6.0cm, width=0.44\textwidth]{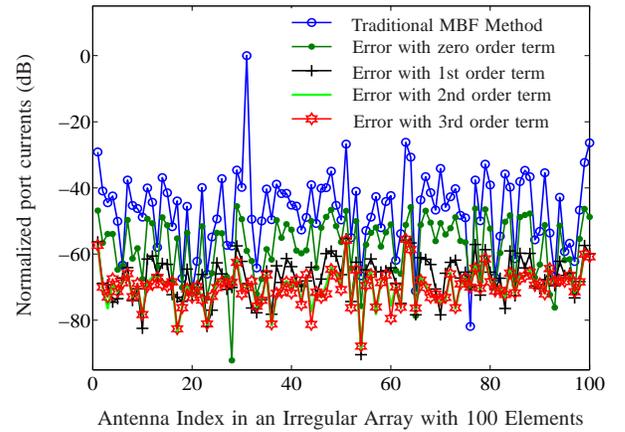}
\end{center}
\caption{Port currents in an irregular array of 100 elements when a single element (indicated by an arrow in the right-down corner of Fig. \ref{fig:array_irrg}) is excited.}
    \label{fig:port_irrg}
\end{figure}
\subsubsection{Port Currents and Radiation Patterns}
The results in terms of normalized port currents obtained with the C-FFT based method are compared in Fig. \ref{fig:port_irrg} with the ones obtained with the traditional MBF method. The port currents are taken when a single element (indicated by an arrow in the right-down corner of Fig. \ref{fig:array_irrg}) is excited and all others are terminated with a $50$ $\Omega$ load. As can be seen from the results, the error level in the computation is below $-30$ dB with second or third order terms. With reference to the interactions in Fig. \ref{fig:evolution_error}, the accuracy in the solution for the case of moderate size arrays can be good enough considering up even down to first order term only. However, as the size of the array increases, one may have to consider upto second or third order terms to have better accuracy. 

\begin{figure}[!t]
\psfragscanon
\begin{center}
\psfrag{a}[][][0.7]{\hspace{15.5em}{Pattern with Traditional MBF Method}}
\psfrag{b}[][][0.7]{\hspace{16.5em}{Pattern with C-FFT Based MBF Method}}
\psfrag{y}[][][0.8]{Directivity (dBi)}
\psfrag{x}[][][0.8]{Angle from Broadside (deg)}
 \includegraphics[height=6.0cm, width=0.44\textwidth]{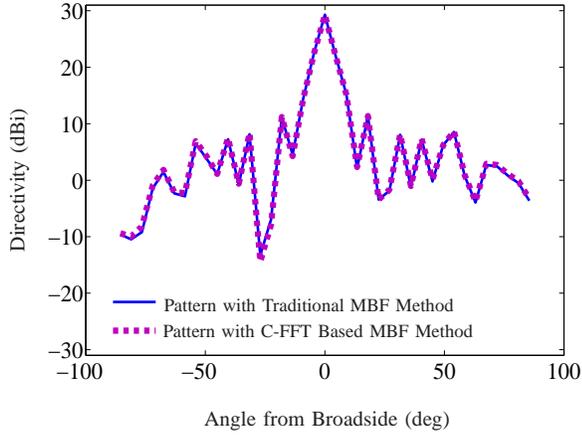}
\end{center}
\caption{An E-plane active array pattern in an irregular array of 100 elements when all the elements are excited.}
\label{fig:pattern_irrg_Eplane}  
\end{figure}

\begin{figure}[!t]
\begin{center}
\psfrag{a}[][][0.7]{\hspace{15.5em}{Pattern with Traditional MBF Method}}
\psfrag{b}[][][0.7]{\hspace{16.5em}{Pattern with C-FFT Based MBF Method}}
\psfrag{y}[][][0.8]{Directivity (dBi)}
\psfrag{x}[][][0.8]{Angle from Broadside (deg)}
\includegraphics[height=6.0cm, width=0.44\textwidth]{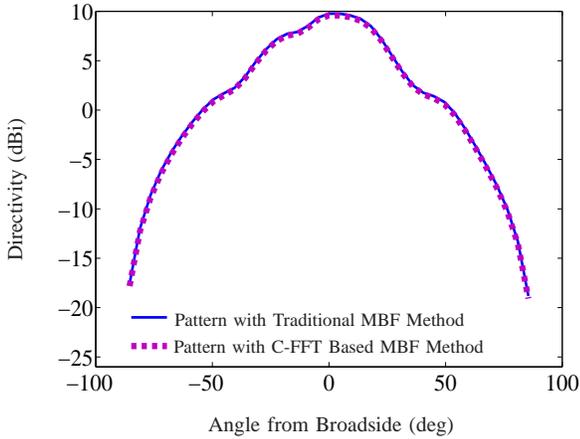}
\end{center}
\caption{An H-plane embedded element pattern when a single element in an array (indicated by an arrow in the middle of Fig. \ref{fig:array_irrg}) is excited.}
\label{fig:pattern_irrg_Hplane}   
\end{figure}

The results for the case of E-plane active array pattern and H-plane embedded element patterns are given in Fig. \ref{fig:pattern_irrg_Eplane} and Fig. \ref{fig:pattern_irrg_Hplane}, respectively. To obtain the active array patterns, all the elements are excited. The embedded element pattern refers to the case when a single element (indicated by an arrow in the middle of Fig. \ref{fig:array_irrg}) is excited and all others are terminated with a $50$ $\Omega$ load. The error level in computation of the directivity values are also very small and lies normally in the range of 0.01 dB to 0.1 dB and between 0.5 dB till 1 dB for angles very far from broadside.

\begin{figure}[!t]
\psfragscanon
\begin{center}
\psfrag{a}[][][0.7]{\hspace{10.75em}{Traditional MBF Method}}
\psfrag{b}[][][0.7]{\hspace{11.75em}{C-FFT based MBF Method}}
\psfrag{c}[][][0.7]{\hspace{9.5em}{Error in the solutions}}
\psfrag{y}[][][0.8]{Normalized port currents (dB)}
\psfrag{x}[][][0.8]{Antenna Index in an Regular Array of $25$ $\times$ $25$ elements}
 \includegraphics[height=5.8cm, width=0.49\textwidth]{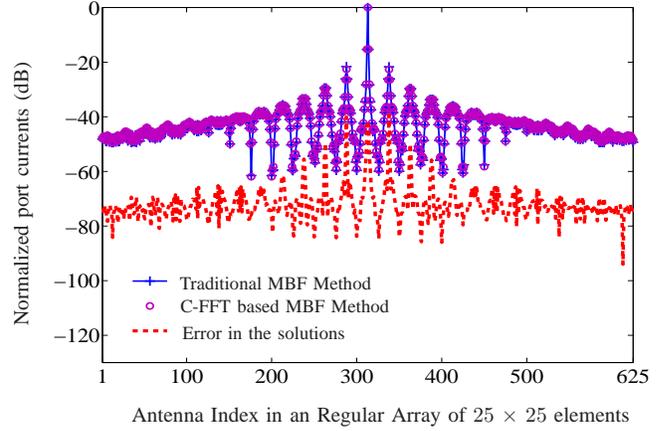}
\end{center}
\caption{Port currents in a large regular $25$ $\times$ $25$ array when a middle element in the array is excited and rest are terminated with a load.}
    \label{fig:port_middle_regular}    
\end{figure}
\subsection{Analysis of a Regular Array of Microstrip Antenna}
\label{sec:regarras}
For this case, a large regular array of size $25$ $\times$ $25$ is considered. The array spacings in X and Y directions are taken equal to $0.58\lambda_0$. The total number of unknowns to solve for is $151875$. The  port currents obtained by using the C-FFT based method are compared with the ones obtained with the traditional MBF method in Fig. \ref{fig:port_middle_regular}. The port currents are taken when a middle element is excited and all others are terminated with a $50$ $\Omega$ load. The comparisons are of excellent accuracy with an error level well below $-30$ dB. It is to be noted here that the accuracy for computations of other quantities like port currents for different excitation conditions, active array pattern and embedded element patterns are of similar accuracy. 
\subsection{Computational Complexity of the Method}
\label{sec:comp}
In this section, we provide the computational complexity of the proposed C-FFT based method for the computation of the substrate-related interactions between MBFs. The computational complexity of each FFT operation is $O(Nlog_2N$). The number $N=\rho^2$ of points in the FFT is delineated as follows. If we denote by $d_{max}$ the maximum size of the array in a given direction in terms of wavelengths, the expression for $\rho$ is arrived at by imposing the Nyquist sampling of $\exp(-ik_{max}x)$ over a distance $d_{max}$ along x:
\begin{equation}
\label{eq:Nterm}
\rho=2\frac{d_{max}k_{max}}{2\pi}
\end{equation} 

Knowing from the Section \ref{sec:FFT-spec} that $k_{max}$ is in the order of $5 k_0 \sqrt{\epsilon_r}$ and defining $S$$>$$1$ as an over-sampling factor, we obtain:
\begin{equation}
\label{eq:Nterm}
\rho=10S\sqrt{\epsilon_r}\frac{d_{max}}{\lambda_0}
\end{equation} 

In our case, we have used an FFT of size $2048$ $\times$ $2048$ including the zero-padding and an over-sampling factor of $S=3$. If one uses $N_t$ Taylor's series terms during FFT applications, the complexity will scale, for $M$ MBFs, as $O(N_tM^2 Nlog_2N)$ for the reduced matrix preparation time for the full arrays. For an example, for a given MBF pair, the total time needed to compute FFTs considering upto 3rd order Taylor's series terms is $1.56$ seconds using the same computing resources as outlined in Section \ref{sec:Results}. Considering $9$ MBFs, leading to $81$ MBFs pair, it takes $126.36$ seconds for all the FFT operations. This corresponds to a fixed amount of tabulation time and does not need to be recomputed for other array configurations as long as the simulation parameters remain the same. There is some additional time needed to fill the reduced matrix using the interpolations. This additional time depends on the complexity of the interpolation routine.  For instance, it takes $1.62$ seconds to compute the $390000$ interpolations using the implemented second-order routine. The computational complexity of $O(N_tM^2 Nlog_2N)$ indicates that the reduced matrix filling time using the C-FFT based MBF method remains independent from the number of elementary basis functions used to represent the antennas. 

The proposed method can be implemented with RAM-based memory, and it is the FFT operations which require most of the memory. For instance, for an FFT of size $2048$ $\times$ $2048$, a RAM memory of $64$ Megabytes is needed. The memory requirement for an FFT operation scale as $N$. Hence, the total memory requirement for full FFT operations considering $N_t$ Taylor's series terms scales as $O(N_t N)$. In addition, some fixed amount of RAM memory is needed for the computations of the Fourier transform of the Green's functions and Macro Basis Functions. For our case, to analyze a large regular arrays of $25$ $\times$ $25$, the total memory requirement was smaller than $1.5$ Gigabytes.

The related computation times are reported in Table \ref{tab:computation}. The values between the parentheses are related to the regular arrays of size $25$ $\times$ $25$. In the case of regular arrays using the traditional MBF method, we make use of the position-related symmetry when filling the reduced matrix. However, to benchmark our fast method, we explicitly computed all the interactions without exploiting the benefits of the position related symmetry. For the case of irregular arrays made of $100$ elements, one needs to compute $10000$ MBF interactions. Out of those, $9900$ are computed with the C-FFT based fast method. The other $100$ interactions are related to the self-impedance reduced matrix, which are computed with the traditional MBF technique. As the elements are identical, this has to be done just once and can be reused requiring just a few seconds. Similarly, for the case of regular arrays of size $25$ $\times$ $25$, one needs to compute $390625$ MBF interactions. Out of those, $390000$ are computed using the C-FFT based method, and the remaining $625$ self-impedance terms using the traditional MBF solution as before. Regarding the multipole calculations, the filling time for the $25$ $\times$ $25$ array is actually over-estimated in the C-FFT column ($4.97$ mins) because$-$contrary to the traditional MBF case in regular arrays and as announced above$-$translational symmetry has not be exploited. The provided time would hence be the same of an
irregular array made of 625 elements. The time needed to obtain the interactions between the MBFs is very small; this allows one to analyze very large arrays in a small amount of time without loosing accuracy in the results.
\begin{table}[htbp]\footnotesize
\caption{Comparison of the computation times for $100$-elements irregular array and for $25$$\times$$25$ regular array (between parentheses), in minutes.}
\label{tab:computation}
\centering
\begin{tabular}{|c||c|c|c|c|c|}
\hline 
\bf{Operations} &\multicolumn{5}{c|}{\bf{Methods/Time (mins)}} \\\cline{2-2}
\cline{3-3}\cline{4-4}\cline{5-5}\cline{6-6}
& \bf{Traditional MBF} & \bf{C-FFT based Method}\\ \hline
MBFs Generation & 1.26 (0.74) & 1.26 (0.74)\\ \hline
FT Tabulation of MBFs & - & 0.51 (0.51) \\ \hline
FFT Time & - & 2.11 (2.11)\\ \hline
\bf{Total Preparation Time} & \bf{1.26 (0.74)} & \bf{3.88 (3.36)} \\ \hline
Reduced Matrix Filling & 743.3 (240.29) & 0.17 (2.19) \\ \hline
MBF-Multipole Part & - & 0.71 (4.97) \\ \hline
Solving Time & 0.001 (0.18) & 0.001 (0.18) \\ \hline
\bf{Total Analysis time} & \bf{744.56 (241.21)} & \bf{4.76 (10.7)}\\ \hline
\end{tabular}
\end{table}

\section{Conclusion}
\label{sec:conclusion}
A Contour-FFT (or C-FFT) based fast MBF method has been proposed to analyse large regular and irregular arrays of antennas printed on a substrate backed by a ground plane. We exploited the benefits of the filtering capabilities of the Macro Basis Functions in spectral domain, which allowed us to limit the range of spectral integration. As discussed in the introduction, the surface wave poles in the layered medium Green's function have been a major accuracy problem and have been tackled previously either by increasing the number of integration points near the poles or by using poles extraction techniques, which may be difficult to control. We hence opted for a classical contour deformation along the radial spectral coordinate. Then, we brought the integration back to rectangular coordinates by proposing a suitable analytical extension of transverse wavenumbers, to be able to exploit the speed-up factor of the FFT algorithm. This transformation in itself did not allow us to apply the FFT directly because of the real exponential factor arising from the contour deformation. We proposed a Taylor's series expansion of this extra factor and moved outside the integral all factors that do not depend on transverse wavenumbers. Then each term becomes suitable for FFT and we demonstrated that just by considering up to $2$nd or $3$rd order terms in the expansion, one could achieve very accurate results in the final solution. The results from the FFT operations provide tables for the interactions between MBFs for any relative positions. Those tables do not need to be recomputed when the array configuration changes.

The C-FFT based method is not only efficient but it also is very simple to implement and has an excellent comparison with the traditional MBF method. The C-FFT based method has a $Nlog_2N$$-$type of complexity in
reduced matrix filling time for the substrate-related terms (the homogeneous-medium part being already efficiently treated in other publications). Using the FFT algorithm on MBFs interactions, a large saving in computation time is achieved. The accuracy obtained in the computations of the port currents and radiation patterns can fulfill the stringent engineering requests. Moreover, the proposed method holds great promise for array optmization problems where one needs an ultra-fast solution at every iteration of the optimization routine while including the mutual coupling effect accurately. The proposed method can be readily extended to FFT-based fast iterative techniques to deal with arbitrarily printed surfaces.

\appendices
\section{Change of Variables used for FFT application}
For maximum spatial limits $\pm(x_m, y_m)$, maximum spectral limits $\pm(k_{xm},k_{ym})$, space-domain indices $k_1$ and $k_2$, and spectral-domain indices $n_1$ and $n_2$, we have, using the notation from Section \ref{sec:fftmethod},
\begin{equation}
K'=(-1)^{(k_1+k_2-2)} \frac{1}{4\pi^2} \Delta k_x \Delta k_y 
\end{equation} 
\begin{equation}
\tilde{I}_{Mf}=(-1)^{(n_1+n_2-2)} \tilde{I}_{M}
\end{equation} 
\begin{equation}
k_x=-k_{xm}+(n_1-1)\frac{\pi}{x_m}
\end{equation} 
\begin{equation}
k_y=-k_{ym}+(n_2-1)\frac{\pi}{y_m}
\end{equation} 
\begin{equation}
\Delta x=-x_m+(k_1-1)\frac{\pi}{k_{x_m}}
\end{equation} 
\begin{equation}
\Delta y=-y_m+(k_2-1)\frac{\pi}{k_{y_m}}
\end{equation} 

\section*{Acknowledgment}
The authors would like to thank Mr. Denis Tihon for his help in Matlab indexing in interpolaion routine and Mr. Abhishek Dangol for the meshing of the Microstrip antenna used as an example.


\ifCLASSOPTIONcaptionsoff
  \newpage
\fi



%
%
%

\bibliographystyle{IEEEtran}
\bibliography{Main_tap_fftFASTSpectral}
\end{document}